\documentclass[12pt,a4paper]{article}

\usepackage{amsthm,amsmath,natbib,amssymb,amsfonts,bm}
\usepackage{epsfig}
\usepackage{placeins}
\usepackage{epsf}
\RequirePackage[dvips]{hyperref}

\def\hang{\hangindent\parindent}
\def\rf{\par\noindent\hang}

\newtheorem{theorem}{Theorem}

\newtheorem*{theorem*}{Theorem}

\theoremstyle{definition}

\theoremstyle{remark}

\begin{document}

\baselineskip=21pt

\begin{center}
 {\bf \Large Confidence intervals in regression centred on the SCAD estimator}
\end{center}

\bigskip

\begin{center}
{\bf \large Davide Farchione and Paul Kabaila$^{\textstyle{^*}}$}
\end{center}

\medskip

\noindent{\it $\textstyle{^*}$Department of Mathematics and
Statistics, La Trobe University, Victoria 3086, \newline Australia}

\bigskip
\noindent{\bf Abstract}
\medskip

Consider a linear regression model. Fan and Li (2001) describe the smoothly clipped absolute deviation (SCAD)
point estimator of the regression parameter vector. To gain insight into
the properties of this estimator, they consider an orthonormal design
matrix and focus on the estimation of a specified component of this vector.
They show that the SCAD
point estimator has three attractive properties.
We answer the question: To what extent can an {\sl interval} estimator, centred on the
SCAD estimator, have similar attractive properties?

\bigskip
\bigskip


\vbox{\vskip 6.5cm}


\noindent $^*$ Corresponding author. Address: Department of
Mathematics and Statistics, La Trobe University, Victoria 3086,
Australia; Tel.: +61-3-9479-2594; fax: +61-3-9479-2466. \newline
{\it E-mail address:} P.Kabaila@latrobe.edu.au.

\newpage

\noindent {\bf 1. Introduction}

\medskip

Consider the linear regression model $Y = X \beta + \varepsilon$, where
$Y$ is a random $n$-vector of responses, $X$ is a known $n \times p$
design matrix with linearly independent columns, $\beta$ is an
unknown $p$-vector and $\varepsilon \sim N(0, \sigma^2 I_n)$, where
$\sigma^2$ is an unknown positive parameter.
In a widely-cited paper, Fan and Li (2001) describe a point estimator
of $\beta$ that they call the smoothly clipped absolute deviation
(SCAD) estimator. The SCAD point estimator
is designed to perform especially well when most of the components of
$\beta$ are believed to be zero (a sparsity type of assumption).

In section 2 of Fan and Li (2001),
to gain insight into the properties of this point estimator, the authors
focus on the estimation of $\beta_i$ (where $i$ is specified) for the case
that the columns of $X$ are orthonormal (cf section 2.2 of Tibshrani, 1996 and
P\"otscher and Schneider, 2010).
This is the scenario that we consider throughout the present paper.
Let sign$(x)$ be equal to $-1$ for $x<0$, 0 for $x=0$ and 1 for $x>0$ and let
$x_+ = \max\{x,0\}$. Let $\hat \beta_i$ denote the least squares estimator of $\beta_i$.
Also let $\hat \Sigma^2$ denote the usual unbiased estimator of $\sigma^2$.
The SCAD estimator of $\beta_i$ is
\begin{equation*}
\tilde{\beta}_i =
\begin{cases}
\text{sign}(\hat \beta_i) \, \big(|\hat \beta_i| - \lambda \big)_+
&\text{if }\ |\hat \beta_i| \le 2 \lambda \\
\big((a-1) \hat \beta_i - \text{sign}(\hat \beta_i) a \lambda \big)/(a-2)
&\text{if }\  2 \lambda < |\hat \beta_i| \le a \lambda \\
\hat \beta_i
&\text{if }\  |\hat \beta_i| > a \lambda
\end{cases}
\end{equation*}
We adopt the proposal of Fan and Li (2001) that $a=3.7$.
For
the purpose of gaining insight into the properties of the SCAD estimator,
these authors suppose that (a) $\sigma$ is known to be 1 and $\lambda$ is a specified
fixed value when they consider the mean square error (m.s.e.) of this estimator and
(b) $\lambda = \lambda_n$ is a non-random sequence that depends on $n$ when they consider
what they call the oracle property. In the present paper, for
the purpose of gaining insight into the properties of {\sl confidence intervals} centred on
the SCAD estimator,
we suppose that $\lambda = \hat{\Sigma} \, \eta$ where $\eta$ is a specified
positive number.

To assess the SCAD point estimator, assume (for the moment) that $\sigma$ is known
and that $\lambda = \sigma \, \eta$, where $\eta$ is a specified
positive number.
We assess the SCAD point estimator
by the ratio (m.s.e. of the SCAD estimator)/(m.s.e. of
least squares estimator), which we call the scaled m.s.e..
This point estimator has the following attractive properties:

\begin{enumerate}

\item[(P1)] It is a continuous function of the data.

\item[(P2)] The scaled m.s.e. converges to 1 as $|\beta_i/\sigma| \rightarrow \infty$.

\item[(P3)] The scaled m.s.e. is substantially less than 1 when
$\beta_i = 0$.

\end{enumerate}

Now consider {\sl interval} estimation of $\beta_i$.
We assess a $1-\alpha$ confidence interval $J$ for $\beta_i$ by the ratio
$E(\text{length of } J)/E(\text{length of usual $1-\alpha$ confidence interval})$,
which we call the scaled expected length.
The corresponding attractive properties for a $1-\alpha$ confidence interval
for $\beta_i$ are the following:

\begin{enumerate}

\item[(I1)] The endpoints of $J$ are continuous functions of the data.

\item[(I2)] The scaled expected length converges to 1 as $|\beta_i/\sigma| \rightarrow \infty$.

\item[(I3)] The scaled expected length is substantially less than 1 when
$\beta_i = 0$.

\end{enumerate}

\noindent Farchione and Kabaila (2008) have already found $1-\alpha$ confidence intervals
that possess all of these attractive properties. These intervals also have the appealing property
that the maximum of the scaled expected length is not too large.
The centres of these interval estimators
do not resemble a SCAD estimator. This suggests that a $1-\alpha$ confidence interval centred on the
SCAD estimator will not be able to have all of the attractive properties (I1), (I2) and (I3).

The SCAD estimator of $\beta_i$ reverts to the least squares estimator $\hat{\beta}_i$ when
$|\hat{\beta}_i| > a \,\hat{\Sigma} \, \eta$. We consider a $1-\alpha$ confidence interval (for $\beta_i$)
centred on this SCAD estimator that, similarly, reverts to the usual $1-\alpha$ confidence interval
for $\beta_i$ when $|\hat{\beta}_i| > a \,\hat{\Sigma} \, \eta$.
This confidence interval has the attractive property (I2).
We will also construct this confidence interval to have the attractive property (I1).
We ask the following question.
To what extent can this confidence interval, centred on the
SCAD estimator, have the property (I3)?
Let $m=n-p$. In Section 3, we consider $1-\alpha = 0.95$  and
the cases (a) $m=200$ (moderately large $m$) and $\eta = 0.5, 1, 2$
and (b) $m=3$ (small $m$) and $\eta = 0.5, 1, 2$. In each of these cases, we show numerically
that this confidence interval, centred on the
SCAD estimator, \textbf{cannot} have the property (I3).
This suggests that this confidence interval cannot have this property more
generally.

The SCAD point estimator may be viewed as being obtained from $\hat \beta_i$, by
a modification determined by $|\hat \beta_i|/\hat \Sigma$. Such a modification seems reasonable
because $|\hat \beta_i|/\hat \Sigma$ may be viewed as a test statistic for testing the null
hypothesis $\beta_i=0$ against the alternative hypothesis $\beta_i \ne 0$.
In the present paper, we consider interval estimators centred at this SCAD estimator, with
width $2 \hat \Sigma \, s(|\hat \beta_i|/\hat \Sigma)$, where the function $s$ is quite
flexible (the constraints on this function are specified in the next section).
This width may be viewed as a modification of a given (non-random) multiple of $\hat \Sigma$, by
a modification determined by $|\hat \beta_i|/\hat \Sigma$.
We use a finite-sample analysis of this
confidence interval; we do not use any asymptotic approximations. To assume that $\sigma^2$ is
known is effectively equivalent to assuming that $n-p$ is large; we do not assume that $\sigma^2$ is
known. We require only that $n-p \ge 1$.
In related work, P\"otscher and Schneider (2010) consider confidence intervals that include in
their interior the hard-thresholding,
LASSO (or soft thresholding) and adaptive LASSO estimators. However, these intervals are constrained
to have a width that is a given (non-random) multiple of $\hat \Sigma$ (or $\sigma$ in the case that they assume
that $\sigma^2$ is known). So, the analysis carried out by P\"otscher and Schneider (2010) is quite
different from the analysis presented in the present paper.

\bigskip

\noindent {\bf 2. The form of the confidence interval centred on the SCAD estimator}

\medskip

Define the quantile $t(m)$ by the requirement that $P(-t(m) \le T \le t(m)) = 1-\alpha$ for
$T \sim t_m$. The usual $1-\alpha$ confidence interval for $\beta_i$ is
\begin{equation*}
I = \big [ \hat \beta_i - t(m) \hat \Sigma, \, \hat \beta_i + t(m) \hat \Sigma \big ].
\end{equation*}
We consider the following confidence interval for $\beta_i$, centred at the SCAD estimator
$\tilde \beta_i$:
\begin{equation*}
J(s) = \big [ \tilde \beta_i - \hat \Sigma \, s(|\hat \beta_i|/\hat \Sigma), \,
\tilde \beta_i + \hat \Sigma \, s(|\hat \beta_i|/\hat \Sigma) \big ],
\end{equation*}
where $s: (0,\infty) \rightarrow (0,\infty)$ is a continuous function that satisfies
$s(x) = t(m)$ for all $x \ge k$, where $k = a \, \eta = 3.7 \, \eta$. This confidence interval has the attractive properties (I1) and (I2).
Farchione and Kabaila (2008) consider $X_1, \ldots, X_n$ independent and identically
$N(\mu, \sigma^2)$ distributed. They consider confidence intervals of the form
\begin{equation*}
\left [ - \tilde{\Sigma} \, c \left(- \frac{\bar{X}}{\tilde{\Sigma}} \right), \,
\tilde{\Sigma} \, c \left( \frac{\bar{X}}{\tilde{\Sigma}} \right) \right],
\end{equation*}
where $\bar{X} = n^{-1} \sum_{i=1}^n X_i$, $\tilde{\Sigma}^2 = (n-1)^{-1} \sum_{i=1}^n (X_i - \bar{X})^2$ and
$c$ is a function satisfying $c(x) \ge -c(-x)$ for all $x \in \mathbb{R}$ (so that the upper
endpoint is always greater than or equal to the lower endpoint). It may be shown that $J(s)$
has a similar form
\begin{equation*}
\left [ - \hat{\Sigma} \, c \left(- \frac{\hat{\beta}_i}{\hat{\Sigma}} \right), \,
\hat{\Sigma} \, c \left( \frac{\hat{\beta}_i}{\hat{\Sigma}} \right) \right].
\end{equation*}
Theorem 1 of Kabaila (2011) implies that if $s$ is chosen such
that $J(s)$ is a $1-\alpha$ confidence interval, with scaled expected length less than 1 when $\beta_i=0$,
then the maximum value of the scaled expected length of $J(s)$ must be greater than 1.

The question that we ask is whether or not we can find a function $s$ such that
$J(s)$ has the property (I3). We do this by minimizing the scaled expected length of $J(s)$
when $\beta_i = 0$, subject to the constraint that the coverage probability of $J(s)$
never falls below $1-\alpha$.

\bigskip

\noindent {\bf 3. Numerical results}

\medskip

\noindent As noted in Appendix A, the scaled expected length and the coverage probability of
$J(s)$ are even functions of $\theta = \beta_i/\sigma$. Let $e(\theta;s)$ denote the scaled
expected length of $J(s)$. To minimize the scaled expected length of $J(s)$
when $\theta=0$ (which is equivalent to $\beta_i = 0$),
subject to the constraint that the coverage probability of $J(s)$
never falls below $1-\alpha$, we use the computationally-convenient expressions described
in Theorem 1 (stated and proved in Appendix A). In Appendix B, we describe briefly how the
coverage probability of $J(s)$ is computed using this theorem.

For computational tractability, we have chosen the function $s$ to be a natural cubic spline with
equally-spaced knots in the interval $[0,k]$ (with a knot at 0 and a knot at $k$).
Remember, $k = a \, \eta = 3.7 \, \eta$. Let these knots
be denoted $x_1, \ldots, x_q$, where $x_1 = 0$ and $x_q = k$. Since we require that $s(x_q) = t(m)$,
the objective function and the constraints for the constrained minimization problem that we consider
are functions of the $q-1$ variables $s(x_1), \ldots, s(x_{q-1})$.

Suppose that $1-\alpha = 0.95$. For $m=200$ (moderately large $m$) and $m=3$ (small $m$)
and for $\eta = 0.5, 1, 2$, we have
computed the function $s$ (specified by $s(x_1), \ldots, s(x_{q-1})$)
that minimizes $e(0;s)$, subject to the constraints that
(a) $s(x)>0$ for all $x \in [0,k]$ and (b)
the coverage probability
of $J(s)$ never falls below $1-\alpha$. Let $s^*$ denote this constrained minimizing value
of the function $s$.  The properties of $s^*$ are summarized in the
Tables 1 and 2 and Figures 1 and 2, below. The function $s^*$ depends on $1-\alpha$, $m$, $k$
and $x_1, \ldots, x_{q-1}$. For notational convenience, this dependence is left implicit.

We implement this coverage constraint in the
computations as follows. It may be shown that, for any reasonable choice of the function $s$,
the coverage probability of $J(s)$ converges to $1-\alpha$ as $\theta \rightarrow \infty$. The
constraints implemented in the computations are that the coverage probability of $J(s)$ is
greater than or equal to $1-\alpha$ for every $\theta$ in a judiciously-chosen finite set of
values. That a given finite set of values of $\theta$ is adequate to the task is judged by checking
numerically,
at the completion of computations, that the coverage probability constraint is satisfied
for all $\theta \ge 0$.

Table 1 presents some properties of this constrained minimizing function $s^*$ for
the case that $m=200$ and $\eta = 0.5, 1, 2$. The number of knots of the cubic spline
$s$ in the interval $[0,k]$ was chosen to be 4, 5 and 6 for each $\eta$. Observe that, for
each value of $\eta$ considered, $e(0;s^*)$ is a decreasing function of the number of knots
and that the decrease from 5 to 6 knots is small.
This table shows that the
confidence interval $J(s)$, which is centred on the SCAD estimator, {\sl cannot} possess the
property (I3) for $m=200$ and these values of $\eta$ and numbers of equally-spaced knots.

\begin{table}[h]\caption{Some properties of the constrained minimizing
function $s^*$ for $m = 200$ and $\eta = 0.5, 1$ and 2.}\label{Tab1}
  \begin{center}
    \begin{tabular}{|c|c|c|c|}
      \multicolumn{4}{c}{$\eta = 0.5$}\\
      \hline
      number of knots & $4$ & $5$ & $6$\\
      \hline
      $e(0;s^*)$ & 1.1609 & 1.1274 & 1.1250\\
      \hline
      $\text{max}_{\theta}e(\theta;s^*)$ & 1.1609 & 1.1274 & 1.1250\\
      \hline
      \multicolumn{4}{c}{}\\
      \multicolumn{4}{c}{$\eta = 1$}\\
      \hline
      number of knots & $4$ & $5$ & $6$\\
      \hline
      $e(0;s^*)$ & 1.2940 & 1.2826 & 1.2825\\
      \hline
      $\text{max}_{\theta}e(\theta;s^*)$ & 1.3936 & 1.3821 & 1.3748\\
      \hline
      \multicolumn{4}{c}{}\\
      \multicolumn{4}{c}{$\eta = 2$}\\
      \hline
      number of knots & $4$ & $5$ & $6$\\
      \hline
      $e(0;s^*)$ & 1.2181 & 1.2155 & 1.2154\\
      \hline
      $\text{max}_{\theta}e(\theta;s^*)$ & 2.1045 & 5.5869 & 5.5272\\
      \hline
    \end{tabular}
  \end{center}
\end{table}

\FloatBarrier

Table 2 presents some properties of the constrained minimizing function $s^*$ for
the case that $m=3$ and $\eta = 0.5, 1, 2$. The number of knots of the cubic spline
$s$ in the interval $[0,k]$ was chosen to be 4, 5 and 6 for each $\eta$.
Observe that, for
each value of $\eta$ considered, $e(0;s^*)$ is a decreasing function of the number of knots
and that the decrease from 5 to 6 knots is small.
This table shows that the
confidence interval $J(s)$, which is centred on the SCAD estimator, {\sl cannot} possess the
property (I3) for $m=3$ and these values of $\eta$ and numbers of equally-spaced knots.

\begin{table}[h]\caption{Some properties of the constrained minimizing
function $s^*$ for $m = 3$ and $\eta = 0.5, 1$ and 2.}\label{Tab2}
  \begin{center}
    \begin{tabular}{|c|c|c|c|}
      \multicolumn{4}{c}{$\eta = 0.5$}\\
      \hline
      number of knots & $4$ & $5$ & $6$\\
      \hline
      $e(0;s^*)$ & 1.0526 & 1.0519 & 1.0511\\
      \hline
      $\text{max}_{\theta}e(\theta;s^*)$ & 1.0759 & 1.0782 & 1.0796\\
      \hline
      \multicolumn{4}{c}{}\\
      \multicolumn{4}{c}{$\eta = 1$}\\
      \hline
      number of knots & $4$ & $5$ & $6$\\
      \hline
      $e(0;s^*)$ & 1.0977 & 1.0966 & 1.0950\\
      \hline
      $\text{max}_{\theta}e(\theta;s^*)$ & 1.3216 & 1.3385 & 1.3464\\
      \hline
      \multicolumn{4}{c}{}\\
      \multicolumn{4}{c}{$\eta = 2$}\\
      \hline
      number of knots & $4$ & $5$ & $6$\\
      \hline
      $e(0;s^*)$ & 1.0824 & 1.0815 & 1.0788\\
      \hline
      $\text{max}_{\theta}e(\theta;s^*)$ & 2.0858 & 2.1650 & 2.1193\\
      \hline
    \end{tabular}
  \end{center}
\end{table}

\FloatBarrier


We now examine the properties of the constrained minimizing function $s^*$ in more detail for
the case that $m=200$, $\eta = 1$ and the cubic spline $s$ has $6$ equally-spaced knots in the interval $[0,k]$.
The top panel of Figure 1 is a plot of the scaled expected length $e(\theta;s^*)$ as a function of $\theta$.
This plot illustrates the fact that every confidence interval of the form $J(s)$ possesses
the attractive property (I2). The bottom panel of this figure is a plot of the coverage probability of
$J(s^*)$ as a function of $\theta$. It is notable that this coverage probability is far above 0.95 for
$\theta \in [0,1]$. We would like to be able to choose the function $s$ so as to
``trade'' this high coverage probability
for a small scaled expected length at $\theta=0$. Evidently,
using a confidence interval of the form $J(s)$, centred on the SCAD estimator, does not allow
this ``trade'' to occur. This is in sharp contrast to the confidence interval of Farchione and Kabaila (2008),
which has coverage probability equal to $1-\alpha$ throughout the parameter space. It appears that by allowing
their confidence interval to have a flexible centre, Farchione and Kabaila (2008) have allowed this ``trade'' to
occur, resulting in a confidence interval that possesses all of the attractive properties (I1), (I2), (I3) and
maximum scaled expected length that is not too large.
Figure 2 is a plot of the constrained minimizing function $s^*$ for this case. The knots of the
cubic spline are denoted by small
circles.

\FloatBarrier

\begin{figure}[h]\hspace{10mm}
    \scalebox{0.6}{\includegraphics[]{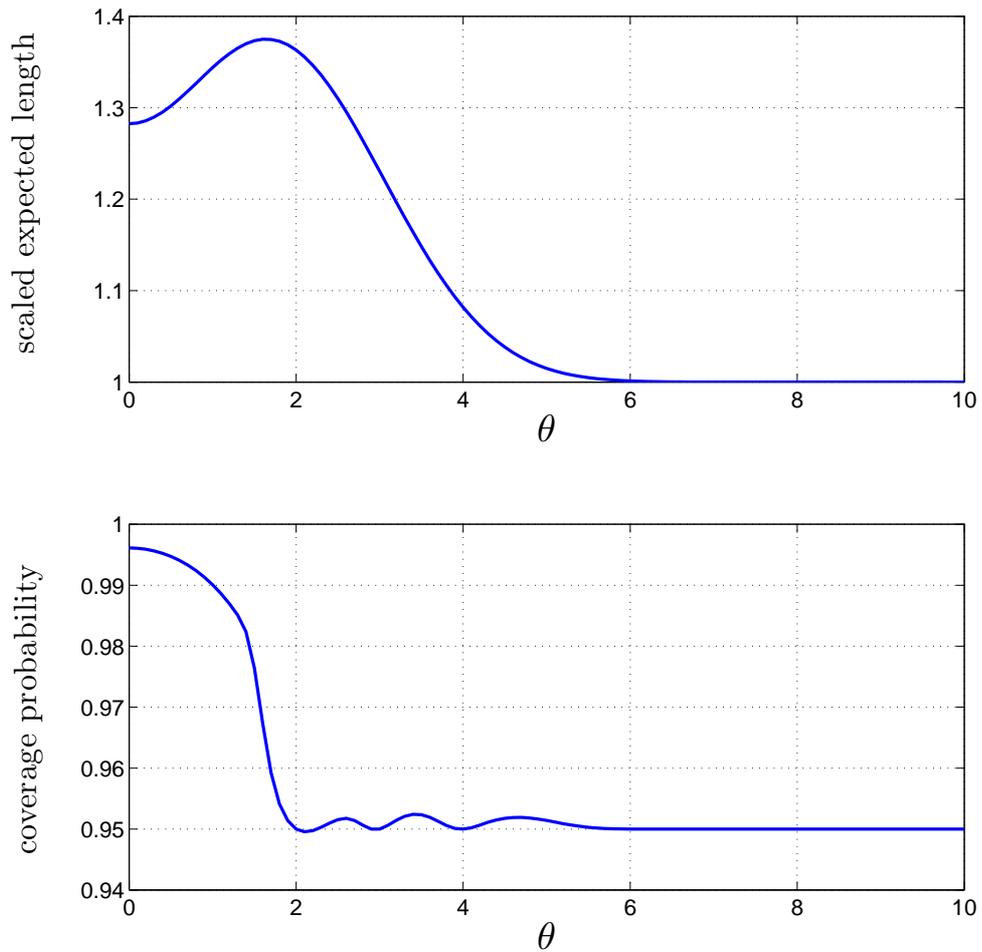}}
    \caption{Properties of the constrained minimizing function $s^*$
    for $m=200$, $\eta = 1$ and the cubic spline $s$ having $6$ equally-spaced knots in the interval $[0,k]$.
    The top panel is a plot of the scaled expected length $e(\theta;s^*)$ as a function of $\theta$.
    The bottom panel is a plot of the coverage probability of
    $J(s^*)$ as a function of $\theta$.}
    \label{Fig1}
\end{figure}

\FloatBarrier

\bigskip

\begin{figure}[h]\hspace{5mm}
    \scalebox{0.6}{\includegraphics[]{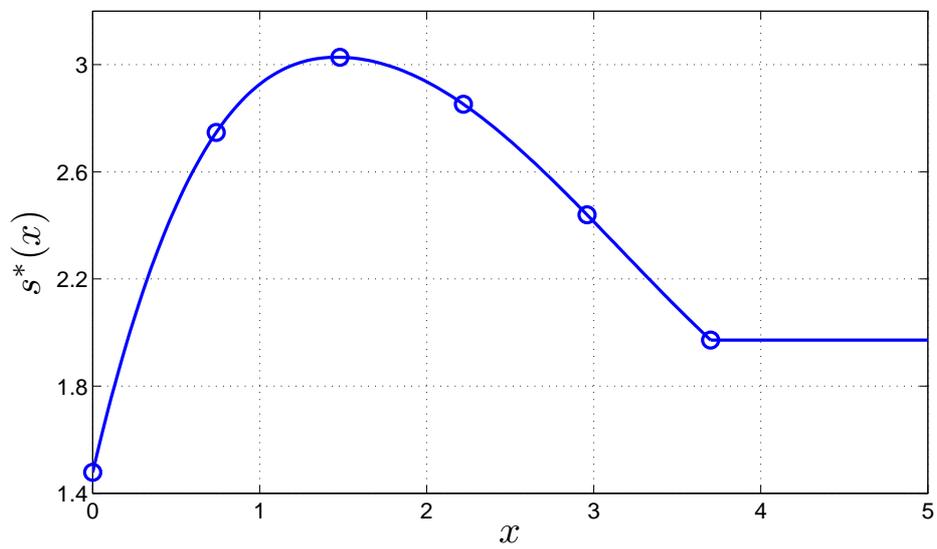}}
    \caption{Plot of the constrained minimizing function $s^*$ for the case that $m=200$,
    $\eta = 1$ and the cubic spline $s$ has $6$ knots in the interval $[0,k]$.}
    \label{Fig2}
\end{figure}

\FloatBarrier

\noindent {\bf 4. Discussion}

\medskip

To gain insight into
the properties of the SCAD estimator, Fan and Li (2001) consider an orthonormal design
matrix and focus on the estimation of a specified component of the regression parameter vector.
We do the same, to gain insight into the properties of confidence intervals
centred on the SCAD estimator. We consider $1-\alpha$ confidence intervals centred on this SCAD estimator
that revert to the usual $1-\alpha$ confidence interval
for the same data values as the SCAD estimator reverts to the least squares estimator.
Our numerical results strongly suggest that these confidence intervals
cannot be constructed to have the the important property (I3).
By contrast, in the context of a multivariate normal mean,
the positive-part James-Stein point estimator
dominates the usual estimator of this mean (using a sum of squared errors loss function).
As shown by Casella and Hwang (1983),
a sphere (with data-dependent radius) centred on this point estimator can be constructed
so as to dominate the usual confidence set for this mean.

\bigskip

\noindent \textbf{Appendix A: Computationally convenient expressions for the scaled expected length
and the coverage probability of $\boldsymbol{J(s)}$}

\medskip

Define $W = \hat{\Sigma}/\sigma$ and let $f_W$ denote the probability density function
of $W$. Let $\theta = \beta_i/\sigma$. Define
\begin{equation*}
h(x) =
\begin{cases}
\text{sign}(x) \, \big(|x| - \eta \big)_+
&\text{if }\ |x| \le 2 \eta \\
\big((a-1) x - \text{sign}(x) a \eta \big)/(a-2)
&\text{if }\  2 \eta < |x| \le a \eta \\
x
&\text{if }\  |x| > a \eta.
\end{cases}
\end{equation*}
We use the notation
\begin{equation*}
{\cal I}({\cal A}) =
\begin{cases}
1 &\text{if } {\cal A} \ \ \text{is true} \\
0 &\text{if } {\cal A} \ \ \text{is false}
\end{cases}
\end{equation*}
where ${\cal A}$ is an arbitrary statement. This is similar to the Iverson bracket notation
(Knuth, 1992). Computationally convenient expressions for the scaled expected
length and coverage probability of $J(s)$ are provided by the following result.

\begin{theorem}

\noindent (a) The scaled expected length of $J(s)$ is equal to
\begin{equation}
\label{sel}
1 + \frac{1}{t(m) E(W)} \int_{-k}^k \big (s(|x|) - t(m) \big) \int_0^{\infty} \phi(wx-\theta) \, w^2 \, f_W(w) \,
dw \, dx.
\end{equation}
For given function $s$, the scaled expected length of $J(s)$ is an even function of $\theta$.

\noindent (b) The scaled expected length of $J(s)$ evaluated at $\theta=0$ is
\begin{equation}
\label{sel_theta_0_comp}
1 + \frac{\sqrt{2/\pi}}{t(m) E(W)} \int_{0}^k \big (s(x) - t(m) \big) \, \left(\frac{m}{x^2+m} \right)^{(m/2)+1} \, dx.
\end{equation}

\noindent (c) Define
\begin{equation*}
b(w;m,k,\theta) =
\begin{cases}
0 \qquad \qquad \text{if } \max(-t(m)w, -kw-\theta) \ge \min(t(m)w, kw-\theta) \\
\Phi \big(\min(t(m)w, kw-\theta) \big) - \Phi \big(\max(-t(m)w, -kw-\theta)\big) \quad \text{otherwise}
\end{cases}
\end{equation*}
where $\Phi$ denotes the $N(0,1)$ cumulative distribution function. The coverage probability of $J(s)$ is equal
to
\begin{align}
\notag
&\int_{-k}^k \int_0^{\infty} {\cal I} \left(h(x) - s(|x|) \le \frac{\theta}{w} \le
h(x) + s(|x|) \right) \, \phi(wx - \theta) \, w \, f_W(w) \, dw \, dx \\
\label{cp}
&+ 1 - \alpha - \int_0^{\infty} b(w;m,d,\theta) \, f_W(w) \, dw
\end{align}
where $\phi$ denotes the $N(0,1)$ probability density function. For given functions $h$ and $s$,
this coverage probability is an even function of $\theta$.

\end{theorem}

\medskip

\noindent \textbf{Proof of part (a)}

\noindent The scaled expected length of $J(s)$ is defined to be
\begin{equation*}
\frac{\text{expected length of } J(s)}{\text{expected length of } I}.
\end{equation*}
This is equal to
\begin{equation}
\label{sel_basic}
\frac{E \big ( s \big(|\hat{\Theta}|/W \big) \, W \big )}{t(m) E(W)}
\end{equation}
where $\hat{\Theta} = \hat{\beta}_i / \sigma$. It follows from Theorem 1(b) of
Kabaila and Giri (2009) that \eqref{sel_basic} is equal to \eqref{sel}.

\medskip

\noindent \textbf{Proof of part (b)}

\noindent It follows from \eqref{sel} that the scaled expected length of $J(s)$ evaluated
at $\theta = 0$ is
\begin{equation}
\label{sel_theta_0}
1 + \frac{1}{t(m) E(W)} \int_{-k}^k \big (s(|x|) - t(m) \big) \int_0^{\infty} \phi(wx) \, w^2 \, f_W(w) \,
dw \, dx.
\end{equation}
Now
\begin{equation*}
\int_0^{\infty} \phi(wx) \, w^2 \, f_W(w) \, dw
= \frac{2 m^{m/2}}{\sqrt{2 \pi}} \frac{1}{\Gamma(m/2) 2^{m/2}}
\int_0^{\infty} w^{m+1} \exp \left ( - \frac{1}{2} \big(m+x^2 \big) w^2 \right ) \, dw
\end{equation*}
where $\Gamma$ denotes the gamma function.
By (A2.1.3) on p.144 of Box and Tiao (1973), this is equal to
\begin{equation*}
\frac{1}{\sqrt{2 \pi}} \left ( \frac{m}{x^2+m} \right )^{(m/2)+1}.
\end{equation*}
\eqref{sel_theta_0_comp} follows from this and \eqref{sel_theta_0}.

\medskip

\noindent \textbf{Proof of part (c)}

\noindent The coverage probability of $J(s)$ is equal to
\begin{equation}
\label{cp_basic}
P \Big ( \tilde{\beta}_i - \hat{\Sigma} \, s \big ( |\hat{\beta}_i|/\hat{\Sigma} \big ) \le \beta_i
\le \tilde{\beta}_i + \hat{\Sigma} \, s \big ( |\hat{\beta}_i|/\hat{\Sigma} \big ) \Big ).
\end{equation}
By the law of total probability, this is equal to
\begin{align*}
&P \Big ( \tilde{\beta}_i - \hat{\Sigma} \, s \big ( |\hat{\beta}_i|/\hat{\Sigma} \big ) \le \beta_i
\le \tilde{\beta}_i + \hat{\Sigma} \, s \big ( |\hat{\beta}_i|/\hat{\Sigma} \big ), |\hat{\beta}_i| \le k \, \hat{\Sigma} \Big ) \\
&+
P \Big ( \tilde{\beta}_i - \hat{\Sigma} \, s \big ( |\hat{\beta}_i|/\hat{\Sigma} \big ) \le \beta_i
\le \tilde{\beta}_i + \hat{\Sigma} \, s \big ( |\hat{\beta}_i|/\hat{\Sigma} \big ), |\hat{\beta}_i| > k \, \hat{\Sigma} \Big ).
\end{align*}
The second term in this sum is equal to
\begin{equation*}
P \Big ( \hat{\beta}_i - t(m) \, \hat{\Sigma} \le \beta_i
\le \hat{\beta}_i + t(m) \, \hat{\Sigma}, |\hat{\beta}_i| > k \, \hat{\Sigma} \Big ).
\end{equation*}
By the law of total probability, this is equal to
\begin{equation*}
1 - \alpha - P \Big ( \hat{\beta}_i - t(m) \, \hat{\Sigma} \le \beta_i
\le \hat{\beta}_i + t(m) \, \hat{\Sigma}, |\hat{\beta}_i| \le k \, \hat{\Sigma} \Big ).
\end{equation*}
Thus \eqref{cp_basic} is equal to
\begin{align*}
&1 - \alpha
+ P \Big ( \tilde{\beta}_i - \hat{\Sigma} \, s \big ( |\hat{\beta}_i|/\hat{\Sigma} \big ) \le \beta_i
\le \tilde{\beta}_i + \hat{\Sigma} \, s \big ( |\hat{\beta}_i|/\hat{\Sigma} \big ), |\hat{\beta}_i| \le k \, \hat{\Sigma} \Big ) \\
&- P \Big ( \hat{\beta}_i - t(m) \, \hat{\Sigma} \le \beta_i
\le \hat{\beta}_i + t(m) \, \hat{\Sigma}, |\hat{\beta}_i| \le k \, \hat{\Sigma} \Big ).
\end{align*}
This is equal to
\begin{align*}
&1 - \alpha
+ P \Big ( \tilde{\Theta} - W \, s \big ( |\hat{\Theta}|/W \big ) \le \theta
\le \tilde{\Theta} + W \, s \big ( |\hat{\Theta}|/W \big ), |\hat{\Theta}| \le k \, W \Big ) \\
&- P \Big ( \hat{\Theta} - t(m) \, W \le \theta
\le \hat{\Theta} + t(m) \, W, |\hat{\Theta}| \le k \, W \Big ).
\end{align*}
where $\hat{\Theta} = \hat{\beta}_i/\sigma$, $\tilde{\Theta} = \tilde{\beta}_i/\sigma$, $\theta = \beta_i / \sigma$ and
$W = \hat{\Sigma}/\sigma$. It may be shown that
\begin{equation*}
\tilde{\Theta} =
\begin{cases}
\text{sign}(\hat{\Theta}) \, \big(|\hat{\Theta}| - W \, \eta \big)_+
&\text{if }\ |\hat{\Theta}| \le 2 W \eta \\
\big((a-1) \hat{\Theta} - \text{sign}(\hat{\Theta}) a W \eta \big)/(a-2)
&\text{if }\  2 \eta < |\hat{\Theta}| \le a \eta \\
\hat{\Theta}
&\text{if }\  |\hat{\Theta}| > a W \eta.
\end{cases}
\end{equation*}
Now define the function $g$ by $\tilde{\Theta} = g(\hat{\Theta}, W)$. Thus
\begin{align}
\notag
&P \Big ( \tilde{\Theta} - W \, s \big ( |\hat{\Theta}|/W \big ) \le \theta
\le \tilde{\Theta} + W \, s \big ( |\hat{\Theta}|/W \big ), \, |\hat{\Theta}| \le k \, W \Big ) \\
\notag
&= E \Big ( {\cal I} \big (g(\hat{\Theta},W) - W \, s(|\hat{\Theta}|/W) \le \theta \le
g(\hat{\Theta},W) + W \, s(|\hat{\Theta}|/W) \big ) \, {\cal I} \big( |\hat{\Theta}| \le k W \big) \Big ) \\
\label{term_cp}
&= \int_0^{\infty} \int_{-kw}^{kw} {\cal I} \big (g(x,w) - w \, s(|x|/w) \le \theta \le
g(x,w) - w \, s(|x|/w) \big ) \, \phi(x-\theta) \, f_W(w) \, dx \, dw.
\end{align}
Now change the variable of integration of the inner integral to $y=x/w$. Thus \eqref{term_cp}
is equal to
\begin{equation}
\label{term_cp_next}
\int_0^{\infty} \int_{-k}^{k} {\cal I} \big (g(wy,w) - w \, s(|y|) \le \theta \le
g(wy,w) - w \, s(|y|) \big ) \, \phi(wy-\theta) \, w \, f_W(w) \, dy \, dw.
\end{equation}
It may be shown that $g(wy,w) = w h(y)$. Thus \eqref{term_cp_next} is equal to
\begin{align}
\notag
& \int_0^{\infty} \int_{-k}^k {\cal I} \left(h(x) - s(|x|) \le \frac{\theta}{w} \le
h(x) + s(|x|) \right) \, \phi(wx - \theta) \, w \, f_W(w) \, dx \, dw  \\
\label{term1_cp_next_next}
&=\int_{-k}^k \int_0^{\infty} {\cal I} \left(h(x) - s(|x|) \le \frac{\theta}{w} \le
h(x) + s(|x|) \right) \, \phi(wx - \theta) \, w \, f_W(w) \, dw \, dx .
\end{align}
Now
\begin{align}
\notag
&P \Big ( \hat{\Theta} - t(m) \, W \le \theta
\le \hat{\Theta} + t(m) \, W, |\hat{\Theta}| \le d \, W \Big ) \\
\label{term2_cp}
&= P \Big (  - t(m) \, W \le Z
\le t(m) \, W, |Z + \theta| \le d \, W \Big )
\end{align}
where $Z = \hat{\Theta} - \theta$, so that $Z \sim N(0,1)$. Observe that \eqref{term2_cp}
is equal to
\begin{equation*}
\int_0^{\infty} P \big( Z \in [-t(m)w,t(m)w] \cap [-dw-\theta, dw-\theta] \big ) \, f_W(w) \, dw.
\end{equation*}
It may be shown that this is equal to
\begin{equation*}
\int_0^{\infty} b(w;m,d,\theta) \, f_W(w) \, dw.
\end{equation*}
Thus the coverage probability is equal to \eqref{cp}. Now \eqref{term1_cp_next_next}
and \eqref{term2_cp} may be shown to be even functions of $\theta$. It follows that
the coverage probability is an even function of $\theta$.

\bigskip

\noindent \textbf{Appendix B: Computation of the coverage probability}

\medskip

By Theorem 1(c), for given functions $h$ and $s$, the coverage probability of $J(s)$ is an
even function of $\theta$. Consequently, we only need to compute this coverage probability
for $\theta \ge 0$. To compute this coverage probability using \eqref{cp}, we need to
compute
\begin{equation}
\label{inner_int_cov}
\int_0^{\infty} {\cal I} \left(h(x) - s(|x|) \le \frac{\theta}{w} \le
h(x) + s(|x|) \right) \, \phi(wx - \theta) \, w \, f_W(w) \, dw
\end{equation}
for given $x \in [-d,d]$. We consider the following 2 cases.

\medskip

\noindent \textbf{Case 1:} $\boldsymbol{\theta = 0}$

\medskip

\noindent In this case,
\begin{align*}
{\cal I} \left(h(x) - s(|x|) \le \frac{\theta}{w} \le
h(x) + s(|x|) \right)
=
\begin{cases}
1 &\text{if } h(x) - s(|x|) \le 0 \text{ and } h(x) + s(|x|) \ge 0 \\
0 &\text{otherwise}.
\end{cases}
\end{align*}
Thus
\begin{equation*}
\eqref{inner_int_cov}
=
\begin{cases}
\int_0^{\infty} \phi(wx) \, w \, f_W(w) \, dw &\text{if } h(x) - s(|x|) \le 0 \text{ and } h(x) + s(|x|) \ge 0 \\
0 &\text{otherwise}.
\end{cases}
\end{equation*}
We find a convenient expression for
\begin{equation}
\label{1D_simple_int}
\int_0^{\infty} \phi(wx) \, w \, f_W(w) \, dw
\end{equation}
as follows. Substituting the formulae for $\phi$ and $f_W$ into \eqref{1D_simple_int},
we find that
\begin{equation}
\label{1D_simpler_int}
\eqref{1D_simple_int} =
\frac{1}{\sqrt{2 \pi}} \frac{2 \, m^{m/2}}{\Gamma(m/2) \, 2^{m/2}} \int_0^{\infty} w^m \exp \left ( - \frac{1}{2}(x^2+m) w^2 \right ) \, dw.
\end{equation}
By (A2.1.3) on p.144 of Box and Tiao (1973),
\begin{equation*}
\eqref{1D_simpler_int} = \frac{1}{\sqrt{\pi}} \frac{\Gamma((m+1)/2)}{\Gamma(m/2)}
\left ( \frac{m}{x^2+m} \right)^{m/2} \frac{1}{\sqrt{x^2+m}}.
\end{equation*}

\medskip

\noindent \textbf{Case 2:} $\boldsymbol{\theta > 0}$

\medskip

\noindent \textbf{Subcase (a):} $\boldsymbol{h(x) - s(|x|) > 0 \textbf{ and } h(x) + s(|x|) > 0}$

\medskip

\noindent In this subcase,
\begin{equation*}
{\cal I} \left(h(x) - s(|x|) \le \frac{\theta}{w} \le
h(x) + s(|x|) \right)
=
{\cal I} \left(\frac{\theta}{h(x) + s(|x|)} \le w \le \frac{\theta}{h(x) - s(|x|)} \right).
\end{equation*}
Thus, in this subcase,
\begin{equation*}
\eqref{inner_int_cov}
=
\int_{\theta/(h(x) + s(|x|))}^{\theta/(h(x) - s(|x|))} \phi(wx-\theta) \, w \, f_W(w) \, dw.
\end{equation*}

\medskip

\noindent \textbf{Subcase (b):} $\boldsymbol{h(x) - s(|x|) \le 0 \textbf{ and } h(x) + s(|x|) > 0}$

\medskip

\noindent In this subcase,
\begin{align*}
{\cal I} \left(h(x) - s(|x|) \le \frac{\theta}{w} \le
h(x) + s(|x|) \right)
&=
{\cal I} \left(0 \le \frac{\theta}{w} \le h(x) + s(|x|) \right) \text{ since } \frac{\theta}{w} > 0 \\
&= {\cal I} \left(\frac{\theta}{h(x) + s(|x|)} \le w < \infty \right).
\end{align*}
Thus, in this subcase,
\begin{equation*}
\eqref{inner_int_cov}
=
\int_{\theta/(h(x) + s(|x|))}^{\infty} \phi(wx-\theta) \, w \, f_W(w) \, dw.
\end{equation*}

\medskip

\noindent \textbf{Subcase (c):} $\boldsymbol{h(x) - s(|x|) < 0 \textbf{ and } h(x) + s(|x|) \le 0}$

\medskip

\noindent In this subcase,
\begin{equation*}
{\cal I} \left(h(x) - s(|x|) \le \frac{\theta}{w} \le
h(x) + s(|x|) \right) = 0 \text{ since } \frac{\theta}{w} > 0.
\end{equation*}
Thus, in this subcase, $\eqref{inner_int_cov}=0$.

\bigskip

\noindent {\bf References}

\smallskip

\rf Box, G.E.P., Tiao, G.C. 1973. Bayesian Inference in Statistical Analysis.
Wiley, New York.

\smallskip

\rf Casella, G., Hwang J.T., 1983. Empirical Bayes confidence
sets for the mean of a multivariate normal distribution. Journal of
the American Statistical Association 78, 688--698.

\smallskip

\rf  Fan, J., Li, R., 2001. Variable selection via nonconcave
penalized likelihood and its oracle properties. Journal of the American
Statistical Association 96, 1348--1360.

\smallskip

\rf   Farchione, D., Kabaila, P., 2008. Confidence intervals for the
normal mean utilizing prior information. Statistics \& Probability Letters
78, 1094--1100.

\smallskip

\rf Kabaila, P., Giri, K., 2009. Confidence intervals in regression
utilizing uncertain prior information. Journal of Statistical Planning and
Inference 139, 3419--3429.

\smallskip

\rf Kabaila, P., 2011. Admissibility of the usual confidence interval
for the normal mean. Statistics \& Probability Letters 81, 352--359.

\smallskip

\rf Knuth, D.E., 1992. Two notes on notation. American Mathematical Monthly  99, 403--422

\smallskip

\rf P\"otscher, B., Schneider, U., 2010. Confidence sets based
on penalized maximum likelihood estimators in Gaussian regression. Electronic
Journal of Statistics 4, 334--360.

\smallskip

\rf Tibshirani, R., 1996. Regression shrinkage and selection via the lasso.
Journal of the Royal Statistical Society, Series B  58, 267--288.

\end{document}